\documentstyle[12pt,bezier]{article}

\title{Global strong solution to a nonlinear  Dirac type equation in one dimension}
\author{Yongqian Zhang
\\ {\small School of Mathematical Sciences, Fudan University,
Shanghai 200433} \\ {\small Key Laboratory of Mathematics for
Nonlinear Sciences }
\\ {\small Email address: yongqianz@fudan.edu.cn} }

\begin{document}
\maketitle

\makeatletter
\renewcommand{\theequation}{\thesection.\arabic{equation}}
\@addtoreset{equation}{section}
\makeatother

\newtheorem{lemma}{Lemma}[section]
\newtheorem{proposition}{Proposition}[section]
\newtheorem{theorem}{Theorem}[section]
\newtheorem{definition}{Definition}[section]
\newtheorem{remark}{Remark}[section]

\begin{abstract}This paper studies a class of nonlinear massless Dirac equations in one dimension, which  include the equations for massless Thirring model and massless Gross-Neveu model. Under the assumptions that the initial data has small charge and is bounded, the global existence  of the strong solution is established. The decay of the local charge is also proved.
\end{abstract}

\section{Introduction}

We consider the nonlinear massless Dirac equations
\begin{equation}\label{eq-dirac}
\left\{ \begin{array}{l} i(u_t+u_x)=G_1(u,v), \\ i(v_t-v_x)=G_2(u,v),
\end{array}
\right.
\end{equation}
with initial data
\begin{equation}\label{eq-dirac-initialv}
(u, v)|_{t=0}=(u_0(x), v_0(x)).
\end{equation}
where $(x,t)\in R^2$, $(u,v)\in \mathbf{C}^2$,  $G_1=\partial_{\overline{u}} W(u,v)+F_1(u,v) $ and $G_2=\partial_{\overline{v}} W(u,v)+F_2(u,v) $  which satisfy the following:
\begin{description}
\item[(A1)] $W(u,v): \mathbf{C}^2\mapsto R$ is a polynomial in variables $|u|^2$ and $|v|^2$;
\item[(A2)] $F_k$ $(k=1,2)$ are polynomials  of cubic order in variables $u,v$ and $\overline{u},\overline{v}$,  moreover there exists a positive constant $c$ such that
\[ |\overline{u}F_1(u,v)|+|\overline{v}F_2(u,v)|\le c |u|^2|v|^2 \]
 for $(u,v)\in \mathbf{C}^2$.
\end{description}
Many physical models verify (A1) and (A2),
for examples,
\begin{itemize}
\item Massless Thirrring model (\cite{thirring})
\[ W=\alpha |u|^2|v|^2, \, \alpha\in R, \, F_1=F_2=0;\]
\item Federbusch model
\[ W=\alpha(|u|^2+|v|^2)|u|^2|v|^2, \, \alpha\in R, \, F_1=F_2=0;\]
\item Massless Gross-Neveu model (\cite{gross-neveu})
\[ (F_1, F_2)=\nabla_{(\overline{u}, \overline{v})} \alpha\big(\overline{u}v +u\overline{v} \big)^2, \, \alpha\in R, \, W=0.\]
\end{itemize}
The nonlinear Dirac equation has been studied by many authors for the massless or massive models, see for examples \cite{bachelot}, \cite{bachelot-2}, \cite{candy}, \cite{chung-p}, \cite{deldado}, \cite{escobedo}, \cite{dias-figueira1}, \cite{dias-figueira2}, \cite{Esteban-Lewin-Sere}, \cite{huh}, \cite{machihara}, \cite{pelimovsky}, \cite{selberg}, \cite{thirring}  and the references therein. The global existence of solution for Cauchy problem has been established in \cite{bachelot}, \cite{deldado}, \cite{dias-figueira1}, \cite{dias-figueira2} in  Sobolev space $H^s$ for suitable $s$. Recently some studies are devoted to the local low regularity solutions in a certain subspace of $C([0,T]; H^s(R))$ for $s>0$; see for instance \cite{machihara}, \cite{pelimovsky} and \cite{selberg} and references therein.In  \cite{huh},  using the explicit representation of the smooth solution by \cite{dias-figueira2} (see also  \cite{machihara}), Huh got the strong solution  in $L^2(R)$ of the nonlinear equations for the massless Thirrring model and Federbusch model. The uniqueness in $L^{\infty}([0,T]: L^2(R)\cap L^4(R))$ and the asymptotic behavior of the field were also established there. In \cite{candy}, Candy used  the strichatz estimates  to prove global existence of $L^2$-strong solution to the nonlinear Dirac equation for massive Thirring model.  For the nonlinear term being a polynomial of quintic or more higher order, the small global solution has been constructed in \cite{pelimovsky} by the end point Strichartz estimates.

Different from the above work, we consider in this paper the case that nonlinear terms include some cubic terms and may have the more complex structure. Such kind of equations include the Massless Gross-Neveu model. To our knowledge, the global existence of strong solution in $L^2$ for Gross-Neveu model is still open \cite{pelimovsky}. Different from the Thirring model,  we have (\ref{eq-dirac1}) for $|u|^2$ and $|v|^2$ instead of the conservation law for the charge $\int_{-\infty}^{\infty} (|u|^2+|v|^2)dx$. The assumption (A2) implies that the equation (\ref{eq-dirac1}) is similar to the approximate conservation law for the interaction between the elementary waves \cite{glimm}; see also \cite{bressan} \cite{dafermos} and \cite{tartar}. Following Glimm's approach, we introduce a Glimm type functional to get the $L^2$ and $L^{\infty}$ bounds of the local $H^1$ solution for  extending the solution globally in time.  Such functional was first given by Glimm \cite{glimm} to establish the global existence of the BV solution for the system of conservation laws, and similar technique has been used to study the general conservation laws, Boltzmann equations and wave maps in $R^{1+1}$, see for instance  \cite{bressan}, \cite{dafermos},  \cite{ha-tzavaras}, \cite{glimm}, \cite{tartar} and \cite{zhou} and the references therein. Our aim is to find the global solutions in the following sense.
\begin{definition}
A function $(u,v)$ is said to be a strong solution to (\ref{eq-dirac}-\ref{eq-dirac-initialv}) on the time interval $[0,T]$ if there exists a sequence of classical solutions $(u^n, v^n)$ to (\ref{eq-dirac}) on $[0,T]\times R$ such that $(u^n,v^n)\to (u,v)$ strongly in $C([0,T]; L^2(R))$ and $(u^n,v^n)|_{t=0}\to (u_0,v_0)$ strongly in $L^2(R)$.
\end{definition}

The main results is stated as follows.
\begin{theorem}\label{thm-L2-existence}
Suppose that $(u_0,v_0)\in L^2(R^1)\cap L^{\infty}(R^1)$ and that the norm $||u_0||_{L^2(R)}+||v_0||_{L^2}$  is small. Then (\ref{eq-dirac}-\ref{eq-dirac-initialv}) admits a unique global strong solution $(u,v)\in C([0,\infty); L^2(R^1))\cap L^{\infty}(R^1\times [0,\infty))$. Moreover, for any constants $a<b$,
\[ \int^b_a (|u(x,t)|^2+|v(x,t)|^2)dx \to 0, \quad t\to +\infty.\]
\end{theorem}

The remaining part of the paper is organized as follows. In section 2, based on the existence of local $H^1-$solutions, we will introduce a Glimm functional for $H^1$ solutions and get the uniform $L^2$ and $L^{\infty}$ bounds of the solution, which are independent of the time and yield the bound of $H^1$ norm for the solution. This enable us to extend the smooth solutions globally in time; see Proposition \ref{prop-global}. In section 3, we take the limit to get the global strong solution in $L^2$ and prove  the decay of the local charge.

\section{Global classical solution }

We consider the case that $(u_0,v_0)\in C_c^{\infty}$ in this section.  (\ref{eq-dirac}) is a semilinear hyperbolic system. Local existence of the solutions in $H^1$ can be proved by the standard methods using the Duhamel's principle and the fixed point argument, see for instance \cite{bachelot-2}, \cite{deldado}  and \cite{tartar}. Now suppose that for $(u_0,v_0)\in C_c^{\infty}$, the initial value problem (\ref{eq-dirac}-\ref{eq-dirac-initialv}) has a smooth solution $(u,v)\in C([0,T_*); H^1(R^1))\cap C^1([0,T_*); L^2(R))$ for some $T_*>0$. Then the support of $(u(t,x),v(t,x))$ is compact for each $t\in [0,T_*)$. To extend the solution across the time $t=T_*$ we only need to show that
 \[ \sup_{0\le t<T_*} \big(||u(t)||_{H^1} +||v(t)||_{H^1}\big)\le C^{\prime}\exp(C^{\prime\prime}T_*) \]
 where the constants $C^{\prime}$ and $C^{\prime\prime}$ are independent of $T_*$.
To this end, we multiply the first equation of (\ref{eq-dirac}) by $\overline{u}$   and the second equation by $\overline{v}$, then
\begin{equation}\label{eq-dirac1}
\left\{\begin{array}{l}
(|u|^2)_t +(|u|^2)_x=2\Re (i\overline{F_1}u), \\  (|v|^2)_t-(|v|^2)_x=2\Re (i\overline{F_2}v).
\end{array}\right.
\end{equation}
Here and in sequel $\Re Z$ denotes the real part of $Z \in \mathbf{C}$.

As in \cite{bony} (see also \cite{ha-tzavaras} and \cite{zhou}), we define the following  functional for the solution for any $t\in [0,T_*)$,
\[ Q(t)=\int\int_{x<y} |u(x,t)|^2|v(y,t)|^2 dxdy,\]
and \[ L(t)=\int^{\infty}_{-\infty} (|u(x,t)|^2+|v(x,t)|^2)dx,\]
\[ D(t)=\int^{\infty}_{-\infty} |u(x,t)|^2|v(x,t)|^2 dx.\]

These functionals were  given in \cite{bony} (see also \cite{ha-tzavaras}). $Q(t) $ is called Bony functional \cite{ha-tzavaras} and is similar to the Glimm interaction potential. The linear combination $L(t)+K Q(t)$ for any constant $K>0$ is called the Glimm functional \cite{glimm}.
Direct computation shows that
\begin{eqnarray*}
\frac{d Q(t)}{dt }  &=& -\int\int_{x<y} (|u(x,t)|^2)_x |v(y,t)|^2 dxdy + \int\int_{x<y} |u(x,t)|^2 (|v(y,t)|^2)_y dxdy \\ & \, & + 2\int\int_{x<y} | \overline{u}F_1(u,v)(x,t)| |v(y,t)|^2 dxdy  \\ &\, &+ 2\int\int_{x<y} |u(x,t)|^2 | \overline{v}F_2(u,v)(y,t)|dxdy  \\
&\le& -2\int^{\infty}_{-\infty} |u(x,t)|^2|v(x,t)|^2 dx  \\ &\, & + 2\int^{\infty}_{-\infty} | \overline{u}F_1(u,v)|dx \int^{\infty}_{-\infty} |v|^2 dy + 2\int^{\infty}_{-\infty} |u|^2 dx\int^{\infty}_{-\infty} | \overline{v}F_2(u,v)|dy \\
\, &\le & \big( -2+2L(t)c \big)\int^{\infty}_{-\infty} |u(x,t)|^2|v(x,t)|^2 dx
\end{eqnarray*}
and
\begin{eqnarray*}
\frac{d L(t)}{dt}\le  2c D(t).
\end{eqnarray*}
Therefore for any constant $K>0$,
\begin{eqnarray}\label{eq-bony-functional-ineq1}
\frac{d L(t)}{dt} +K\frac{d Q(t)}{dt }\le \big\{ (-2+2L(t)c)K+2c \big\}D(t).
\end{eqnarray}

\begin{lemma}\label{lemma-bony-functional} There exists  constants $\delta>0$ and $K>0$ such that for the initial data satisfying $L(0)+KQ(0)\le \delta$ there holds the following
\begin{eqnarray}\label{eq-bony-functional-ineq2}
\frac{d L(t)}{dt} +K\frac{d Q(t)}{dt }\le -D(t)
\end{eqnarray}
for $t\in [0,T_*)$.
Therefore,
\begin{eqnarray}\label{eq-bony-functional-ineq3}
L(t)+K Q(t) + \int_0^T D(\tau) d\tau \le L(0)+K Q(0)
\end{eqnarray}
for $t\in [0,T_*)$.
\end{lemma}
{\it Proof.} First choose constants $\delta$ and $K$ so that
\[ -2+2\delta c<-1/2, \quad -K/2+2c<-2.\]
If $L(0)+KQ(0) \le \delta$, then
$(-2+2L(0)c)K+2c \le -2.$
Let \[ T_0=\sup\{T; (-2+2L(t)c)K+2c<-1, \, t\in[0, T)\}.\]
We want to show that $T_0=T_*$. Otherwise, suppose that $T_0<T_*$, then with (\ref{eq-bony-functional-ineq1}) we have (\ref{eq-bony-functional-ineq3}) for $t\in[0, T_0) $. Therefore,  by the continuity of $L(t)$ and $Q(t)$ for $0\le t<T_*$, we have the following
\[ L(t)+KQ(t)\le \delta, \quad  (-2+2L(t)c)K+2c \le -2, \quad t\in [0, T_0],\]
which leads to the contradiction to the definition of $T_0$ by using again the continuity of $L(t)$ for $0\le t<T_*$.  Thus the proof is complete. $\Box$

\begin{figure}[h]
\begin{center}
\unitlength=1mm
\begin{picture}(130,75)
\put(40,40){$\Omega(x_0,t_0)$}
\put(10,43){$\Gamma_L(x_0,t_0)$}
\put(60,45){$\Gamma_R(x_0,t_0)$}
\put(85,20){\line(-1,1){40}}
\put(5,20){\line(1,1){40}} \put(43,63){$(x_0,t_0)$}
\put(0,20){\vector(1,0){100}}\put(102,20){$x$}
\put(0,16){$(x_0-t_0, 0)$ } \put(78,16){$(x_0+t_0,0)$}
\end{picture}
\caption{Domain $\Omega(x_0,t_0)$}\label{fig-domain}
\end{center}
\end{figure}
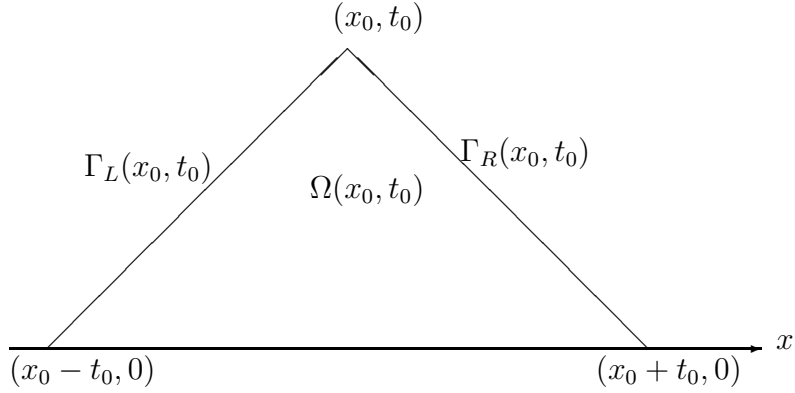

Next we shall use the above estimates to get the $L^{\infty}$ bound of the solution for $0<t<T_*$. Denote \[\Omega(x_0,t_0)=\{(x,t): 0<t<t_0, \, x_0-(t_0-t)<x<x_0+(t_0-t)\},\] and
\[ \Gamma_R(x_0,x_0)=\{(x,t): 0<t<t_0, \, x=x_0+(t_0-t)\}, \]
\[ \Gamma_L(x_0,x_0)=\{(x,t): 0<t<t_0, \, x=x_0-(t_0-t)\}, \]
see Fig. \ref{fig-domain}.

\begin{lemma}\label{lemma-integral}
For any $t_0\in (0,T_*)$, there hold
\begin{eqnarray*}
\int_{\Gamma_R(x_0,t_0)} |u|^2 \le \sqrt{2}c \int\int_{\Omega(x_0,t_0)}|u|^2|v|^2 dxdt +\frac{1}{\sqrt{2}} \int^{x_0+t_0-t}_{x_0-t_0+t} |u_0(x)|^2 dx
\end{eqnarray*}
and
\begin{eqnarray*}
\int_{\Gamma_L(x_0,t_0)} |v|^2 \le \sqrt{2}c \int\int_{\Omega(x_0,t_0)}|u|^2|v|^2 dxdt +\frac{1}{\sqrt{2}} \int^{x_0+t_0-t}_{x_0-t_0+t} |v_0(x)|^2 dx.
\end{eqnarray*}
\end{lemma}
{\it Proof.} Integrating (\ref{eq-dirac1}) over $\Omega(x_0,t_0)$, then by Green formula we have
\begin{eqnarray*}
\sqrt{2}\int_{\Gamma_R(x_0,t_0)} |u|^2 &= &\int\int_{\Omega(x_0,t_0)}2\Re (i\overline{F_1}u)dxdt + \int^{x_0+t_0-t}_{x_0-t_0+t} |u_0(x)|^2 dx\\
 \, &\le & 2c \int\int_{\Omega(x_0,t_0)}|u|^2|v|^2 dxdt + \int^{x_0+t_0-t}_{x_0-t_0+t} |u_0(x)|^2 dx
\end{eqnarray*}
and
\begin{eqnarray*}
\sqrt{2}\int_{\Gamma_L(x_0,t_0)} |v|^2 &= &\int\int_{\Omega(x_0,t_0)}2\Re (i\overline{F_2}v)dxdt + \int^{x_0+t_0-t}_{x_0-t_0+t} |v_0(x)|^2 dx\\
 \, &\le & 2c \int\int_{\Omega(x_0,t_0)}|u|^2|v|^2 dxdt + \int^{x_0+t_0-t}_{x_0-t_0+t} |v_0(x)|^2 dx.
\end{eqnarray*}
The proof is complete. $\Box$

\begin{lemma}\label{lemma-bound}
Let $\delta$ and $K$ be the constants given in Lemma \ref{lemma-bony-functional}. Then for any constants $a$ and $b$ with $a<b$ there hold the following
\[ \sup_{0\le t<T_*, a<x<b} \big( |u(x,t)|^2+|v(x,t)|^2 \big)\le C_1\sup_{ a-T_*<x<b+T_*} \big( |u_0(x)|^2+|v_0(x)|^2 \big), \]
and
\[ \int_a^b \big( |u(x,t)|^2+|v(x,t)|^2 \big) dx \le C_1\int_a^b \big(|u_0(x-t)|^2+|v_0(x+t)|^2\big) dx,\]
where $C_1=\exp((2\sqrt{2}c^2+\sqrt{2}c)\delta)$.
\end{lemma}
{\it Proof.}
 Using the characteristic method, by (\ref{eq-dirac1}) and (A2) we have
\begin{eqnarray*}
\frac{d}{dt}|u(x_0+(t-t_0),t)|^2 & =& 2\Re (i\overline{F_1}u) \\ &\le& 2c |u(x_0+(t-t_0),t)|^2|v(x_0+(t-t_0),t)|^2
\end{eqnarray*}
and
\begin{eqnarray*}
\frac{d}{dt}|v(x_0-(t-t_0),t)|^2 & =& 2\Re (i\overline{F_2}v) \\ &\le& 2c |u(x_0-(t-t_0),t)|^2|v(x_0-(t-t_0),t)|^2,
\end{eqnarray*}
for any $(x_0,t_0)$ with $a-(T_*-t_0)<x_0<b+(T_*-t_0)$.
Then,
\begin{eqnarray*}
|u(x_0,t_0)|^2 &\le & |u(x_0-t_0,0)|^2\exp\big( 2c\int_0^{t_0}|v(x_0+(\tau-t_0),\tau)|^2 d\tau\big)\\
 &=& |u(x_0-t_0,0)|^2\exp\big( 2c\int_{\Gamma_L(x_0,t_0)} |v|^2\big) \\
 &\le & |u_0(x_0-t_0)|^2\exp\big( I_1(x_0,t_0)\big)
\end{eqnarray*}
and
\begin{eqnarray*}
|v(x_0,t_0)|^2 &\le & |v(x_0+t_0,0)|^2\exp\big(2c \int_0^{t_0}|u(x_0-(\tau-t_0),\tau)|^2 d\tau\big)\\
 &=& |v(x_0+t_0,0)|^2\exp\big( 2c\int_{\Gamma_R(x_0,t_0)} |u|^2\big) \\
 &\le & |v(x_0+t_0,0)|^2\exp\big( I_2(x_0,t_0)\big),
\end{eqnarray*}
where the estimates  for $\int_{\Gamma_R(x_0,t_0)} |u|^2 $  and $\int_{\Gamma_L(x_0,t_0)} |v|^2 $ in Lemma \ref{lemma-integral} are used, and
\begin{eqnarray*}
I_1(x_0,t_0)=2\sqrt{2}c^2 \int\int_{\Omega(x_0,t_0)}|u|^2|v|^2 dxdt+\sqrt{2}c\int_{x_0-t_0}^{x_0+t_0} |u_0|^2dx, \\
I_2(x_0,t_0)=2\sqrt{2}c^2 \int\int_{\Omega(x_0,t_0)}|u|^2|v|^2 dxdt+\sqrt{2}c\int_{x_0-t_0}^{x_0+t_0} |v_0|^2dx.
\end{eqnarray*}
 Therefore these lead to the desired result by Lemma \ref{lemma-bony-functional}. The proof is complete. $\Box$

Now, to get the $H^1$ bound of the solution, we differentiate the equations (\ref{eq-dirac}) with respect to $x$, then
\begin{eqnarray*}
\left\{ \begin{array}{l} i(u_{xt}+u_{xx})=\nabla_{u,v,\overline{u},\overline{v}}G_1(u,v)\cdot(u_x,v_x,\overline{u}_x,\overline{v}_x ), \\ i(v_{xt}-v_{xx})=\nabla_{u,v,\overline{u},\overline{v}}G_2(u,v)\cdot(u_x,v_x,\overline{u}_x,\overline{v}_x ),
\end{array}
\right.
\end{eqnarray*}
which lead to the following
\begin{eqnarray*}
\sup_{0\le t<T_*} \big(||u(t)||_{H^1}+ ||v(t)||_{H^1}\big)\le \big(||u_0||_{H^1}+ ||v_0||_{H^1}\big)\exp(C^{\prime}T_*)<\infty,
\end{eqnarray*}
where the constant $C^{\prime}$ depends only on $L^{\infty}$ bounds of $(u,v)$ given in Lemma \ref{lemma-bound}.
Therefore, the local solution could be extended globally. Since $C_c^{\infty}(R^1)$ is dense in $H^1(R^1)$, with the above estimates, we get the following.
\begin{proposition}\label{prop-global}
Suppose that $(u_0,v_0)\in H^1(R^1)$ and that $||u_0||_{L^2(R^1)}+||v_0||_{L^2(R^1)}$ is small enough. Then (\ref{eq-dirac}-\ref{eq-dirac-initialv}) has a unique global classical solution $(u,v)\in C([0,+\infty); H^1(R^1))\cap C^1([0,+\infty); L^2(R^1))\cap L^{\infty}(R^1\times [0,\infty))$. Moreover, there exists some positive constant $C_2$ such that
\begin{eqnarray*}
\sup_{t>0, x\in R^1} \big( |u(x,t)|^2+|v(x,t)|^2 \big)\le C_2\sup_{ x\in R^1} \big( |u_0(x)|^2+|v_0(x)|^2 \big), \\
\int_a^b \big( |u(x,t)|^2+|v(x,t)|^2 \big) dx \le C_2\int_a^b \big(|u_0(x-t)|^2+|v_0(x+t)|^2\big) dx,
\end{eqnarray*}
for any $a<b$.
\end{proposition}

\section{Proof of the main result}

Let $(u_{j,0},v_{j,0})\in C^{\infty}_c(R^1)$ be a sequence of smooth functions such that $(u_{j,0},v_{j,0})\to (u_0,v_0) $ strongly in $L^2$ and such that
 $\sup_j ||(u_{j,0},v_{j,0})||_{L^{\infty}} <\infty$.

Denote $(u_j,v_j)$ be  the $H^1$-solutions of (\ref{eq-dirac}) with initial data $(u_{j,0},v_{j,0})$. Then by Proposition \ref{prop-global} we derive the following
\begin{eqnarray}\label{eq-sequence-bound}
\sup_{t>0, j} (||u_j(t)||_{L^{\infty}} +||v_j(t)||_{L^{\infty}})<\infty.
\end{eqnarray}
Now for any $k,j$
\[
\left\{ \begin{array}{l} (u_k-u_j)_t+(u_k-u_j)_x)=iG_1(u_j,v_j)-iG_1(u_k,v_k), \\ (v_k-v_j)_t-(v_k-v_j)_x)=iG_2(u_j,v_j)-iG_2(u_k,v_k),
\end{array}
\right.
\]
which together with (\ref{eq-sequence-bound}) gives the following by the energy estimate method
\begin{eqnarray*}
&||(u_k-u_j)(t)||_{L^2}^2+||(v_k-v_j)(t)||_{L^2}^2 \, \quad \,  \\  & \le \exp(O(1)t) \big( ||(u_k-u_j)(0)||_{L^2}^2+||(v_k-v_j)(0)||_{L^2}^2\big),
\end{eqnarray*}
where the bounds of $O(1)$ depend only on the $L^{\infty}$ bounds of the solutions $ (u_j,v_j)$ $(j=1,2,\cdots)$ for $t\ge 0$. This yields the convergence of the sequence $\{(u_j,v_j)\}_{j=0}^{\infty}$ to a function $(u,v)$ in $C([0,T); L^2(R^1))$ for any $T>0$. Then $(u,v)$ is a global strong solution to (\ref{eq-dirac}-\ref{eq-dirac-initialv}).
 Moreover, by using Proposition \ref{prop-global} and by repeating the  argument of taking limit, we have
\begin{eqnarray*}
\int_a^b \big( |u(x,t)|^2+|v(x,t)|^2 \big) dx &\le & C_1\int_a^b \big(|u_0(x-t)|^2+|v_0(x+t)|^2\big) dx \\
 & =& C_1 \big( \int_{a+t}^{b+t} |u_0(y)|^2 dy +\int_{a-t}^{b-t} |v_0(y)|^2dy \big).
\end{eqnarray*}
Then the decay of the local charge  follows.

In the same way, we can prove the following for any strong solutions $(u^{\prime},v^{\prime})$ and $(u^{\prime\prime},v^{\prime\prime})$,
\begin{eqnarray*}
||(u^{\prime}-u^{\prime\prime})(t)||_{L^2}^2+||(v^{\prime}-v^{\prime\prime})(t)||_{L^2}^2 \le \exp(O(1)t) \big( ||(u^{\prime}-u^{\prime\prime})(0)||_{L^2}^2+||(v^{\prime}-v^{\prime\prime})(0)||_{L^2}^2\big)
\end{eqnarray*}
for $t\ge 0$,
where the bounds of $O(1)$ depend only on the $L^{\infty}$ bounds of the solutions. This leads to the uniqueness of the strong solution.
The proof of Theorem \ref{thm-L2-existence} is complete.

\section*{ Acknowledgement}

This work was partially  supported by NSFC Project 11031001 and 11121101, by the 111 Project
B08018 and by Ministry of Education of China.


\begin{thebibliography}{9}

\bibitem{bachelot} A. Bachelot, {\it Global existence of large amplitude solutions for nonlinear massless Dirac equation,} Portucaliae Math.{\bf 46} Fasc. Supl. (1989), 455-473.

\bibitem{bachelot-2} A. Bachelot, {\it Global Cauchy problem for semilinear hyperboloc systems with nonlocal interactions. Applications to Dirac equations,} J. Math. Pures. Appl. {\bf 86} (2006), 201-236.

\bibitem{bony} J. M. Bony, {\it Solution globales born\'{e}es pour les mod\`{e}les discrets de l'\'{e}quation de Boltzmann, en dimension 1 d'espace,} Acte des Journ\'{e}es E.D.P. \`{a} Saint-Jean-de-Monts, Ecole Polytechnique (1987), XVI-1-XVI-9.

\bibitem{bressan} A. Bressan, {\it Hyperbolic systems of conservation laws:
 The one-dimensional Cauchy problem,} Oxford University Press
 Inc., New York, 2000.

\bibitem{candy} T. Candy, {\it Global existence for an $L^2$ critical nonlinear Dirac equation in one dimension,}  Adv. Differential Equations {\bf 16} No. 7-8 (2011), 643-666.

\bibitem{chung-p} M. Chungunova and D. Pelinovsky, {\it Block-diagnonalization of the symmetric first-order coupled-mode system,} SIAM J. Appl. Dyn.System  {\bf 5} (2006), 66-83.

\bibitem{dafermos} C. M. Dafermos, {\it Hyperbolic Conservation Laws in Continuum Physics}, Springer-Verlag, Berlin, 2010.

\bibitem{dias-figueira2} J.P. Dias and M. Figueira,{\it Time decay for the solutions of a nonlinear Dirac equation in one space dimension,} Ricerche Mat. {\bf 35} (1986), 309-316.

\bibitem{deldado} V. Delgado, {\it Global solution of the Cauchy problem for the (classical) coupled Maxwell-Dirac and other nonlinear Dirac equations,} Proc. Amer. Math. Soc. {\bf 69}(2)(1978), 289-296.

\bibitem{dias-figueira1} J. P. Dias and M. Figueira, {\it Sur l'existence d'une solution globale pour une equation de Dirac non lin\'{e}aire avec masse nulle,} C.R. Acad. Sci. Paris, S\'{e}rie I, {\bf 305} (1987), 469-472.

\bibitem{dias-figueira2} J. P. Dias and M. Figueira, {\it Remarque sur le probl\`{e}me de Cauchy pour une equation de Dirac non lin\'{e}aire avec masse nulle,} Portucaliae Math. {\bf 45}(4) (1988), 327-335.

\bibitem{escobedo} M. Escobedo and L. Vega, {\it A semilinear Dirac equation in $H^s(R^3)$ for $s>1$,} SIAM J. Math. Anal. {\bf 28}(2) (1997), 338-362.

\bibitem{Esteban-Lewin-Sere} M. J. Esteban, M. Lewin and E. S\'{e}r\'{e}, {\it Variational methods in relativistic quantum mechanics,} Bulletin A.M.S. {\bf 45} No. 4 (2008), 535-598.

\bibitem{glimm} J. Glimm, {\it Solution in the large
for nonlinear systems of conservation laws,} Comm.
Pure Appl. Math. {\bf 18}(1965), 695-715.

\bibitem{gross-neveu} D.J. Gross and A. Neveu, {\it Dynamical symmetry breaking in asymptotically free field theories,} Phys. Rev. D {\bf 10} (1974), 3235-3253.

\bibitem{ha-tzavaras} S-Y Ha and A. E. Tzavaras, {\it Lyapunov functionals and $L^1$-stability for discrete velocity Boltzmann equations,} Commun. Math. Phys. {\bf 239} (2003), 65-92.


\bibitem{huh} H. Huh, {\it Global strong solution to the Thirring model in critical space,} J. Math. Anal. Appl. {\bf 381} (2011), 513-520.

\bibitem{machihara} S. Machihara and T. Omoso, {\it The explicit solutions to the nonlinear Dirac equation and Dirac-Klein-Gordon equation,} Ricerche Mat. {\bf 56} (1)(2007), 19-30.

\bibitem{pelimovsky} D. Pelinovsky, {\it Survey on global existence in the nonlinear Dirac equations in one dimension,} arXiv:1011.5925vl [math-ph] 26 Nov 2010.

\bibitem{selberg} S. Selberg and A. Tesfahun, {\it Low regularity well-posedness for some nonlinear Dirac equations in one space dimension,} Differential Integral Equations {\bf 23}(3-4)(2010), 265-278.

\bibitem{tartar} L. Tartar, {\it From Hyperbolic Systems to Kinetic Theory,} Springer 2008.

\bibitem{thirring} W.E. Thirring, {\it A soluble relativistic field theory,} Ann. Phys. {\bf 3}(1958), 91-112.

\bibitem{zhou} Y. Zhou, {\it Uniqueness of weak solutions in 1+1 dimensional wave maps,} Math. Z. {\bf 232}(1999), 707-719.


\end{thebibliography}
\end{document}